\documentclass{amsart}

\usepackage{pb-diagram}
\usepackage{amsmath}
\usepackage{amssymb}
\usepackage[russian,english]{babel}

\theoremstyle{plain}
\newtheorem{theorem}{Theorem}[section]
\newtheorem{lemma}{Lemma}[section]

\newtheorem{definition}{Definition}[section]
\newtheorem{corollary}{Corollary}[section]

\numberwithin{equation}{section}

\title[On the extension of the Toeplitz algebra by isometries]{On the extension of the Toeplitz algebra by isometries}
\author{T.A. Grigoryan, \ E.V. Lipacheva, \ V.H. Tepoyan}

\begin{document}

\maketitle

\begin{center}
\textit{ Kazan State Power Engineering University }
\end{center}

\textbf{Abstract. }{ We introduce the notion of $\pi$-extension of the semigroup $\mathbb{Z}_+$ and study the extensions of the Toeplitz algebras by
 isometric operators. We show that when the action of the Toeplitz algebra is irreducible all such extensions generate the same algebra,
 i.e. there are no non-trivial extensions of the Toeplitz algebra. Also we provide the examples of the non-trivial extensions of the 
Toeplitz algebra in case its representation is reducible.. }



\section*{Introduction}

We introduce the notion of $\pi$-extension of the semigroup of non-negative integers $\mathbb{Z}_+$ (see definition \ref{D1}). 
We study the properties of $C^*$-algebras generated by $\pi$-extension of the semigroup $\mathbb{Z}_+$. 
Also the concept of the inverse $\pi$-extension of semigroup $\mathbb{Z}_+$ is introduced (see definition \ref{D3}).
 It is proved that if $\pi$ is an irreducible representation of the semigroup $\mathbb{Z}_+$, then there is no non-trivial inverse
 $\pi$-extension. And in the case $\pi$ is  reducible there exists a non-inverse $\pi$-extension.
 On the other hand we show that for any isometric representation of $\mathbb{Z}_+$ there always exists a non-inverse $\pi$-extension.

We also study the extensions of the Toeplitz algebra generated by the inverse $\pi$-extensions of the semigroup $\mathbb{Z}_+$. It is
 shown that among such extensions there are extensions which are the tower of the nested Toeplitz algebras inductive limit of which is the $C^*$-algebra considered in the work of Douglas \cite{Douglas}.

The authors are sincerely thankful to Suren A. Grigoryan for useful discussions and for valuable advice.

\section{Preliminaries}\label{par2}

Let $H^2$ be the Hardy space on the unit circle $S^1$ and $\{e^{in\theta}\}_{n=0}^{\infty}$ be its orthonormal basis. A linear operator
$$T:H^2\rightarrow H^2 \ \ Te^{in\theta}=e^{i(n+1)\theta}$$
is called a shift operator. Closed in the operator norm $C^*$-subalgebra of the algebra of bounded linear operators $B(H^2)$ on $H^2$
generated by  shift operator and its adjoint is the Toeplitz algebra and denoted by $\mathcal{T}$.

This operator generates an isometric representation
$$\pi_0:\mathbb{Z}_+\rightarrow B(H^2)$$
of the semigroup of non-negative integers $\mathbb{Z}_+$: $\pi_0(n)=T^n$.

Throughout the paper by $\pi_0$ we denote this representation of the semigroup $\mathbb{Z}_+$.

Consider the non-unitary isometric representation
$$\pi:\mathbb{Z}_+\rightarrow B(H),$$
where $H$ is a Hilbert space. Let $C^*_{\pi}(\mathbb{Z}_+)$ be a $C^*$-subalgebra of the algebra $B(H)$ generated by the operator $\pi(n)$ and its adjoint $\pi^*(n)$. By Coburn's theorem \cite{Coburn} the algebra $C^*_{\pi}(\mathbb{Z}_+)$ is $\ast$-isomorphic to the Toeplitz algebra $\mathcal{T}$. So in what follows  algebras of the form $C^*_{\pi}(\mathbb{Z}_+)$ are called Toeplitz algebras. Note that if we consider the semigroup $\mathbb{Z}_+\backslash\{1\}$, then there are canonically non-isomorphic representations \cite{Raeburn}.

Let $Is(H)$ be a semigroup of isometric operators in $B(H)$.

\begin{definition}\label{D1}
The subsemigroup $M$ of  $Is(H)$ is called $\pi$-extension of the semigroup $\mathbb{Z}_+$, if
\begin{enumerate}
	\item $\pi(\mathbb{Z}_+)\subset M;$
	\item $\pi(1)T=T\pi(1)$ for all  $T$ in $M$.
\end{enumerate}
\end{definition}

We denote by $C^*_{\pi}(M)$ the $C^*$ algebra generated by  semigroup $M$. This algebra is called 
\emph{$\pi$-extension of the algebra $C^*_{\pi}(\mathbb{Z}_+)$ by the semigroup $M$.}

In this paper we investigate  $C^*$-algebras generated by $\pi$-extensions of the semigroup $\mathbb{Z}_+$.

\section{Inverse Representations}

Suppose that $S$ is an additive commutative cancellative semigroup. We denote by $\Delta_S$ the set of all isometric representations of
  $S$. For $\pi\in\Delta_S$ define the semigroup $S^{\pi}$ generated by finite product of  operators from the semigroup $\pi(S)$ and from its adjoint semigroup $\pi(S)^*$.

If the semigroup $S^{\pi}$ consists of  isometries and partial isometries, then conjugation turns $S^{\pi}$ into an inverse semigroup. Recall that the semigroup in which every element $a$ has a unique element $a^*$ such that
$$aa^*a=a, \ \ a^*aa^*=a^*$$
is called \emph{ inverse semigroup.}

\begin{definition}\label{D2}
The representation ${\pi}\in{\Delta_{S}}$ is called inverse if $S^{\pi}$ is an inverse semigroup.
\end{definition}

\emph{ The regular isometric representation} is the map ${\pi}:S\rightarrow{B(l^2(S))}$, $a\mapsto {\pi}(a)$ defined as follows:

\[
({\pi}(a)f)(b) =
\begin{cases}
f(c), & \text{if $b=a+c$ for $c\in S$;} \\
0, & \text{otherwise.}
\end{cases}
\]

The $C^*$-algebra generated by the regular isometric representation of the semigroup $S$ is called \emph{ reduced semigroup $C^*$-algebra } 
and is denoted by $C^{*}_{red}(S)$ \cite{Jang, Grigoryan, Jang3}.

\begin{theorem}\label{1}
The regular representation of the semigroup $S$ is inverse (see \cite{AT, Grigoryan2}).
\end{theorem}

Let
$$\pi_1:S\rightarrow C^*_{\pi_1}(S)\subset B(H_1) \mbox{ and } \pi_2:S\rightarrow C^*_{\pi_2}(S)\subset B(H_2)$$
be representations of the semigroup $S$. The algebras $C^*_{\pi_1}(S)$ and $C^*_{\pi_2}(S)$ \emph{ are canonically isomorphic} if there exists $\ast$-isomorphism $\tau:C^*_{\pi_1}(S)\rightarrow C^*_{\pi_2}(S)$ such that the diagram
\[ \begin{diagram}
\node{S} \arrow[2]{e,t}{\pi_1}
         \arrow{se,b}{\pi_2}
\node[2]{C^*_{\pi_1}(S)} \arrow{sw,l}{\tau}\\
\node[2]{C^*_{\pi_2}(S)}
\end{diagram}\]
is commutative.

Define a partial order on the semigroup $S$: $a\prec b$, if $b=a+c$. It is said that the semigroup is totally ordered if either $a\prec b$ or $b\prec a$ for any $a\neq b$ from $S$.

\begin{corollary}
Let $S$ be  totally ordered. Then the semigroups $S^{\pi}, \ \pi\in\Delta_S$ are inverse and isomorphic.
\begin{proof}
By Murphy's theorem \cite{Murphy} if $\pi:S\rightarrow B(H_{\pi})$ is a non-unitary isometric representation then the algebra $C^*_{\pi}(S)$ is canonically isomorphic to $C^*_{red}(S)$. Canonical isomorphism of these algebras generates $\ast$-isomorphism of the semigroups $S^{red}$ and $S^{\pi}$, where $S^{red}$ is a semigroup in $C^*_{red}(S)$ generated by the regular representation $red:S\rightarrow B(l^2(S))$. The semigroup $S^{red}$ is inverse so $S^{\pi}$ is an inverse semigroup.
\end{proof}
\end{corollary}

Here is an example of a non-inverse representation.

Let $\pi:\mathbb{Z}_+\rightarrow B(H^2)$ be a representation of the semigroup $\mathbb{Z}_+$ described in section \ref{par2} such that $\pi(n)$ is an  operator of multiplication by a function $e^{in\theta}$.

Every inner function $\Phi(z)$ also defines an isometric multiplicative operator $T_{\Phi}$:
$$T_{\Phi}f=\Phi f.$$

\begin{theorem}\label{2}
Let $\widetilde{\pi}:\mathbb{Z}_+\times\mathbb{Z}_+\rightarrow B(H^2)$ be a representation, which
 maps $(n,0)\mapsto e^{in\theta}$ and $(0,m)\mapsto \Phi^m$, where $\Phi$ is an arbitrary inner function not in
 $\{e^{in\theta}\}_{n=0}^{\infty}$. Then $\widetilde{\pi}$ is a non-inverse representation, i.e. 
$(\mathbb{Z}_+\times\mathbb{Z}_+)^{\widetilde{\pi}}$ is a non-inverse semigroup.

\begin{proof}
An arbitrary inner function has the following form:
$$\Phi(z)=B(z)S(z),$$
where $B(z)$ is the Blaschke product and $S(z)$ is the inner singular function.
In general, $B(z)$ equals zero when $z=0$, so we can say that:
$$\Phi(z)=z^nB_1(z)S(z),$$
where $B_1(z)$ is the Blaschke product which is not equal to 0 at zero. We fix this $n$. Note that $\widetilde{\pi}(k,m)=\pi(k)T_{\Phi}^m$.

For further proof we need to show that
$$\widetilde{\pi}(0,1)\widetilde{\pi}^*(0,1)\widetilde{\pi}(n+1,0)\widetilde{\pi}^*(n+1,0)\neq\widetilde{\pi}(n+1,0)\widetilde{\pi}^*(n+1,0)\widetilde{\pi}(0,1)\widetilde{\pi}^*(0,1) \mbox{ (see \cite{Clifford}) },$$
\begin{equation}\label{ninv}
\mbox{ or } T_{\Phi}T^{*}_{\Phi}{\pi}(n+1){\pi}^{*}(n+1)\neq {\pi}(n+1){\pi}^{*}(n+1)T_{\Phi}T^{*}_{\Phi}.
\end{equation}

Compute  the left and right sides of  inequality (\ref{ninv}) on $e^{in\theta}$. Since ${\pi}^{*}(n+1)(e^{in\theta})=0$ then clearly
$$T_{\Phi}T^{*}_{\Phi}{\pi}(n+1){\pi}^{*}(n+1)(e^{in\theta})=0.$$

Consider $T^{*}_{\Phi}e^{in\theta}$. We show that all  Fourier coefficients of this function in the decomposition
 $\sum_{k=0}^{\infty}{c_ke^{ik\theta}}$ except  zero coefficient are equal to 0. For this we compute scalar product on the elements $e^{ik\theta}$, where $k>0$
$$(e^{ik\theta},T^{*}_{\Phi}e^{in\theta})=(\Phi e^{ik\theta},e^{in\theta})=(e^{ik\theta}B_1(e^{i\theta})S(e^{i\theta}),1)=
\int_{S^{1}}{e^{ik\theta}B_1(e^{i\theta})S(e^{i\theta})d\mu}=0.$$

Now compute the zero coefficient
$$(1,T^{*}_{\Phi}e^{in\theta})=(\Phi,e^{in\theta})=(e^{in\theta}B_1(e^{i\theta})S(e^{i\theta}),e^{in\theta})$$
$$=(B_1(e^{i\theta})S(e^{i\theta}),1)=\int_{S^{1}}{B_1(e^{i\theta})S(e^{i\theta})d\mu}=B_1(0)S(0).$$

Thus $T^{*}_{\Phi}e^{in\theta}=B_1(0)S(0)1$. Hence
$$T_{\Phi}T^{*}_{\Phi}e^{in\theta}=B_1(0)S(0)\Phi=B_1(0)S(0)e^{in\theta}B_1(e^{i\theta})S(e^{i\theta}).$$
Then it is easy to see that:
$${\pi}^{*}(n+1)T_{\Phi}T^{*}_{\Phi}e^{in\theta}={\pi}^{*}(1)B_1(0)S(0)B_1(e^{i\theta})S(e^{i\theta})\neq 0,$$
and consequently
$$({\pi}(n+1){\pi}^{*}(n+1)T_{\Phi}T^{*}_{\Phi})(e^{in\theta})\neq 0.$$

We have used the Cauchy integral formula from the theory of complex functions.

Thus the inequality (\ref{ninv}) is satisfied and therefore $\widetilde{\pi}$ is a non-inverse representation.
\end{proof}

\end{theorem}

\section{$\pi_0$-extension of the Semigroup $\mathbb{Z}_+$}

\begin{lemma}\label{Berl}
Every isometric operator in $\pi$-extension of the semigroup $\mathbb{Z}_+$ can be represented through a single inner function.
\begin{proof}
Suppose that $T$ is an arbitrary isometric operator in $\pi_0$-extension of the semigroup $\mathbb{Z}_+$. We denote:
$$\mathcal{Q}_0=TH^2\subset H^2.$$
Let us show that $\mathcal{Q}_0$ is invariant under translations $\pi_0(1)$. Since $T$ commutes with $\pi_0(1)$  we have
$$\pi_0(1)\mathcal{Q}_0=\pi_0(1)TH^2=T\pi_0(1)H^2\subset TH^2=\mathcal{Q}_0.$$
Then, Beurling's theorem \cite{Garnett}, if $\mathcal{Q}_0\neq{0}$, then there exists an inner function $\Phi$ such that:
$$\mathcal{Q}_0=\Phi H^2,$$
and that $\Phi$ is unique up to a constant factor.

Consider $\Phi^{'}=T\cdot 1$. Since $||(T\cdot 1)h||$ coincides with  $||h||$ for any $h$ in $H^2$, then we get that $\Phi^{'}$ is an inner function, so that:
$$Th=T_{\Phi^{'}}h, \ h\in H^2.$$
Thus, by  Beurling's theorem we obtain that any isometric operator $T$ is represented by the unique inner function.
\end{proof}

\end{lemma}

If $M$ is a $\pi$-extension of the semigroup $\mathbb{Z}_+$, then by lemma \ref{Berl} for each isometric operator $T\in M$ there exists
 a unique inner function $\Phi$, such that the operator $T$ is a multiplicative operator on $\Phi$. Define
$$M^{'}=\{\Phi; \ T_{\Phi}\in M\}.$$

\begin{theorem}
Suppose that $M$ is the $\pi_0$-extension of the semigroup $\mathbb{Z}_+$. Then the following conditions are equivalent:
\begin{enumerate}
	\item $C^{*}_{\pi_0}(\mathbb{Z}_+)=C_{\pi_0}^*(M);$
	\item $M$ is the subsemigroup of the semigroup of finite Blaschke products.
\end{enumerate}

\begin{proof}
Let $M^{'}$ be a subsemigroup of the semigroup of finite Blaschke products. Let's prove that
$C^{*}_{\pi}(\mathbb{Z}_+)=C^*_{\pi}(M)$.

Every finite Blaschke product defines an isometric multiplicative operator and uniformly approximated by finite linear combinations of functions $\{e^{in\theta}\}_{n=0}^{\infty}$. Therefore if $B(z)$ is the finite Blaschke product, then the operator $T_{B}$ belongs to the algebra $C^{*}_{\pi}(\mathbb{Z}_+)$. Consequently,
$$C^{*}_{\pi}(\mathbb{Z}_+)=C^*_{\pi}(M).$$
Now let us prove the converse. Let $T\in M$. Then the operator $T$ corresponds to  some inner function $\Phi\in M^{'}$, such that for any $h\in H^2$:
$$Th=\Phi h.$$

Since $C^{*}_{\pi}(\mathbb{Z}_+)=C^*_{\pi}(M)$, then $T\in C^{*}_{\pi}(\mathbb{Z}_+)$. Then for any $\varepsilon >0$ there exists a finite linear combination of monomials $W_i$, consisting of the degrees of $\pi(1)$ and $\pi^*(1)$, such that
$$||T-\sum_{i}\lambda_iW_i||<\varepsilon.$$
It is easy to check that such monomial can be written as the product of ${\pi}(n){\pi}^*(m)$.

Using $W_i={\pi}(n_i){\pi}^*(m_i)$ in the above inequality, we obtain
$$||T-\sum_{i}\lambda_{i}{\pi}(n_i){\pi}^*(m_i)||<\varepsilon.$$
Let $m_0$ be a maximum of $m_i$. Then
$$||T\pi{(m_0)}-\sum_{i}\lambda_{i}{\pi}(n_i){\pi}^*(m_i)\pi{(m_0)}||<||T-\sum_{i}\lambda_{i}{\pi}(n_i){\pi}^*(m_i)||\cdotp||\pi{(m_0)}||<\varepsilon.$$
Consequently,
$$||T\pi{(m_0)}-\sum_{i}\lambda_{i}{\pi}(n^{'}_i)||<\varepsilon.$$

Because the operator $T$ is a multiplicative operator of multiplication by $\Phi$ and the operator $\pi(n)$ of multiplication by function $e^{in\theta}$, then
$$||\Phi e^{im_0\theta}-\sum_{i}\lambda_{i}e^{in^{'}_i\theta}||<\varepsilon.$$
So, we have
$$||\Phi-\sum_{i}\lambda_{i}e^{i(n^{'}_i-m_0)\theta}||<\varepsilon.$$

Since $\varepsilon$ is arbitrary,  $\Phi$ is uniformly approximated by continuous functions and therefore it is a continuous function on the circle. Since any inner function that is continuous on the circle is a finite Blaschke product, we conclude that $\Phi$ is a finite Blaschke product.
\end{proof}

\end{theorem}

\section{Inverse $\pi$-extension}

We denote by $\mathbb{Z}_+^{\pi}$ the involutive semigroup generated by $\pi(\mathbb{Z}_+)$ and $\pi(\mathbb{Z}_+)^*$. Note that $\mathbb{Z}_+^{\pi}$ is a bicyclic semigroup. All irreducible representations of this semigroup are described in \cite{Arzumanyan}. Let $M$ be a $\pi$-extension of the semigroup $\mathbb{Z}_+$. Denote by $\mathcal{M}^*$ the semigroup generated by $M$ and $M^*$.

\begin{definition}\label{D3}
We call the $\pi$-extension of the semigroup $\mathbb{Z}_+$ inverse, if $\mathcal{M}^*$  is an inverse semigroup.
\end{definition}

Let $\pi:\mathbb{Z}_+\rightarrow B(H^2)$ be the representation of the semigroup $\mathbb{Z}_+$, described in section \ref{par2}. Then we have the following statement.

\begin{theorem}\label{tinv}
$\mathcal{M}^*$ is inverse iff $\mathcal{M}^*=\mathbb{Z}_+^{\pi}$.
\begin{proof}
Let $\mathcal{M}^*=\mathbb{Z}_+^{\pi}$. Every element of the semigroup $\mathbb{Z}_+^{\pi}$ can be represented as the
 product $\pi(n)\pi^*(m)$. Inverse is $\pi(m)\pi^*(n)$. So $\mathbb{Z}_+^{\pi}$ and therefore $\mathcal{M}^*$ is an inverse semigroup.

Conversely, suppose that $\mathcal{M}^*\neq\mathbb{Z}_+^{\pi}$. Then according to Lemma \ref{Berl} every isometric operator $T\in M\backslash\pi(\mathbb{Z}_+)$ 
 is represented through the unique inner function: $T=T_{\Phi}$. Similarly to the theorem \ref{2} it can be shown that in this case $\mathcal{M}^*$ is a non-inverse semigroup.
\end{proof}

\end{theorem}

Consider an arbitrary isometric representation $\pi:\mathbb{Z}_+\rightarrow B(H)$. Denote by $H_0=\ker{\pi^*(1)}$. It is clear that $H_0$ is the Hilbert subspace of $H$.

\begin{theorem}
Let $\pi:\mathbb{Z}_+\rightarrow B(H)$ be an isometric representation of the semigroup $\mathbb{Z}_+$, such that the subspace
$H_0=\ker{\pi^*(1)}$ is not one dimensional. Then there exists an inverse $\pi$-extension $M$ of the semigroup $\mathbb{Z}_+$, such that $\mathbb{Z}_+^{\pi}$ is the proper involutive subsemigroup of the involutive semigroup $\mathcal{M}^*$.
\begin{proof}
Let us show that $ \ pi $ is a reducible isometric representation.
Denote $H_1=\pi(1)H_0, \dots, H_n=\pi(n)H_0, \dots$. The subspaces $H_0, H_1, H_2, \dots$ are pairwise orthogonal. Indeed, we verify that $H_n\bot H_m$ $(m>n)$. Calculate the following
$$(\pi(n)h_0,\pi(m)h_0)=(\pi^*(m-n)h_0,h_0)=0.$$

Thus
$$H=\oplus_{n=0}^{\infty}{H_n},$$
where subspace $H_0$  and, consequently, $H_n, n>0,$ are not one-dimensional.

Let $\{e_0^{(j)}\}, j=1, 2, \dots,$ be an either finite or infinite basis in $H_0$. Suppose that
 $e_n^{(j)}=\pi(n)(e_0^{(j)})$. Then $\{e_n^{(j)}\}, j=1, 2, \dots,$ is a basis $H_n$, and
 $\{e_n^{(j)}\}, j=1, 2, \dots; n=0, 1, 2, \dots$ is a basis in $H$.

Consider the subspace in $H$ with the basis $\{e_n^{(j)}\}_{n=0}^{\infty}$. On this basis, the operator $\pi(1)$ acts as a shift operator
$$\pi(1)(e_n^{(j)})=e_{n+1}^{(j)}.$$

This subspace is the Hardy space, denoted by $H_j^2$.

Thus, the space $H$ is represented as a direct sum of a finite or infinite Hardy subspaces 
$$H=\oplus_{j}{H_j^2},$$
and the representation $\pi$ is irreducible on each of the $H_j^2$.

Since $H_0$ is not one dimensional subspace, then the direct sum $\oplus_{j}{H_j^2}$ contains at least two components. We fix two Hardy spaces $H_1^2$ and $H_2^2$ with basis $\{e_n^{(1)}\}_{n=0}^{\infty}$ and $\{e_n^{(2)}\}_{n=0}^{\infty}$,
respectively.

Consider the operator $T\in B(H)$, which in $H_1^2\oplus H_2^2$ acts as follows
$$T(e_n^{(1)})=e_n^{(2)}; \ T(e_n^{(2)})=e_{n+1}^{(1)}; \ n=0, 1, 2, \dots,$$
and in $H_j^2 \ (j\neq 1, j\neq 2): \ T=\pi(1)$.

Note that on $H_1^2\oplus H_2^2$: $\pi(1)=T^2$.

It follows that $\pi(1)T=T\pi(1)$.

Similarly, it is easy to see that $\pi^*(1)T^*=T^*\pi^*(1)$, but $T^*\pi(1)\neq\pi(1)T^*$, $T\pi^*(1)\neq\pi^*(1)T$.

Let $M$ be the semigroup generated by $\pi(\mathbb{Z}_+)$ and by operator $T$. Obviously it is $\pi$-extension of the semigroup
 $\mathbb{Z}_+$ and $\mathbb{Z}^{\pi}_+\subsetneq\mathcal{M}^*$.

Note that $\mathcal{M}^*$ is an inverse semigroup. This follows from the fact that any monomial in $\mathcal{M}^*$, as can be seen, has the following form

$$T^kT^{*l}\pi(n)\pi^*(m)T^sT^{*p},$$

and inverse to it will be a monomial

$$T^pT^{*s}\pi(m)\pi^*(n)T^lT^{*k}.$$
\end{proof}

\end{theorem}

\section{Inverse Extension of Toeplitz Algebra}

Suppose we have two singular inner functions:
$$\Phi_1=\exp{\frac{e^{i\theta}+1}{e^{i\theta}-1}} \ \mbox{and} \ \Phi_{t}=\exp{t\frac{e^{i\theta}+1}{e^{i\theta}-1}},$$
where $t$ is a positive real number.

Function $\Phi_1$ corresponds to an isometric multiplicative operator $T_{\Phi_1}$, and the function $\Phi_{t}$ to an operator $T_{\Phi_{t}}$.

Note that, isometric operators $T_{\Phi_1}$ and $T_{\Phi_{t}}$, $t>0$, are mapping the Hardy space $H^2$ again in $H^2$, i.e.
$$T_{\Phi_{t}}=PT_{\Phi_{t}} \ \mbox{ and } \ T_{\Phi_1}=PT_{\Phi_1},$$
where $P$ is a projection from $L^2(S^{1}d\mu)$ on $H^2(S^{1}d\mu)=H^2$.

The adjoint operator, for example, to the operator $T_{\Phi_{t}}$ is $T^*_{\Phi_{t}}=PT_{\Phi^{-1}_{t}}=PT_{\Phi_{-t}}$. The operators $T^*_{\Phi_{t}}$ and $T^*_{\Phi_{1}}$ are not isometric.

By the theorem of Coburn \cite{Coburn} the $C^*$-algebra generated by an isometric operator $T_{\Phi_1}$ is canonically isomorphic to 
the Toeplitz algebra. Denote by $\mathcal{T}_1$ the Toeplitz $C^*$-algebra generated by the operator $T_{\Phi_1}$, and by
 $\mathcal{T}_t$ the Toeplitz $C^*$-algebra generated by $T_{\Phi_{t}}$.

Consider the $C^*$-algebra generated by the operators $T_{\Phi_1}$ ш $T_{\Phi_{t}}$. We denote it by $C^{*}(T_{\Phi_1},T_{\Phi_t})$. It is clear that
$$\mathcal{T}_1\subset C^{*}(T_{\Phi_1},T_{\Phi_t}).$$

If $t$ is р positive rational number, then we have the following lemma.

\begin{lemma}\label{tirat}
Let $t=\frac{m}{n}$, where $m$ and $n$ are relatively prime integers. Then $C^{*}(T_{\Phi_1},T_{\Phi_t})\cong\mathcal{T}_{\frac{1}{n}}$, where $\mathcal{T}_{\frac{1}{n}}$ is the Toeplitz $C^*$-algebra generated by the operator $T_{\Phi_{\frac{1}{n}}}$.

\begin{proof}
Clearly the following equalities hold
$$T_{\Phi_{\frac{m}{n}}}=T^m_{\Phi_{\frac{1}{n}}} \ \mbox{ and } \ T_{\Phi_{1}}=T^n_{\Phi_{\frac{1}{n}}}.$$
This means that $C^{*}(T_{\Phi_1},T_{\Phi_{\frac{m}{n}}})\subseteq\mathcal{T}_{\frac{1}{n}}$.

On the other hand, if $m$ and $n$ are relatively prime integers then from the \emph{Euclidean algorithm} follows that there are integers $k$ and $l$, such that
$$nk+ml=1 \mbox{ or } k+\frac{lm}{n}=\frac{1}{n}.$$

There are two cases:
\begin{enumerate}
	\item $k>0, \ l<0$;
	\item $k<0, \ l>0$.
\end{enumerate}

In the first case it is easy to see that:
$$T_{\Phi_{1/n}}=(T^{*}_{\Phi_{m/n}})^{|l|}T^{k}_{\Phi_1},$$
and in the second case:
$$T_{\Phi_{1/n}}=(T^{*}_{\Phi_{1}})^{|k|}T^{l}_{\Phi_{m/n}}.$$
Indeed, for example in the first case we have
$$(T^{*}_{\Phi_{m/n}})^{|l|}T^{k}_{\Phi_1}=(PT_{\Phi_{-m/n}})^{|l|}T^{k}_{\Phi_1}=PT_{\Phi_{-m/n}}PT_{\Phi_{-m/n}}\dots PT_{\Phi_{-m/n}}T^{k}_{\Phi_1}.$$
Whereas
$$T_{\Phi_{-m/n}}T^{k}_{\Phi_1}=T_{\Phi_{-m/n}}T_{\Phi_k}=T_{\Phi_{k-\frac{m}{n}}} \mbox{ and } k-\frac{m}{n}>0,$$
then $PT_{\Phi_{-m/n}}T^{k}_{\Phi_1}=T_{\Phi_{k-\frac{m}{n}}}$.

Similarly,
$$PT_{\Phi_{-m/n}}T_{\Phi_{k-\frac{m}{n}}}=T_{\Phi_{k-\frac{2m}{n}}},$$
etc.. Finally, we obtain
$$(PT_{\Phi_{-m/n}})^{|l|}T^{k}_{\Phi_1}=T_{\Phi_{k-\frac{|l|m}{n}}}=T_{\Phi_{k+\frac{lm}{n}}}=T_{\Phi_{1/n}}.$$

This means that $\mathcal{T}_{\frac{1}{n}}\subseteq C^{*}(T_{\Phi_{1}},T_{\Phi_{\frac{m}{n}}})$.

Thus, we obtain $C^{*}(T_{\Phi_1},T_{\Phi_{\frac{m}{n}}})\cong\mathcal{T}_{\frac{1}{n}}$.

\end{proof}

\end{lemma}

From the lemma~\ref{tirat} it follows that $C^*$-algebra $\mathcal{T}_{\frac{1}{n}}$ is the inverse extension of the $C^*$-algebra $\mathcal{T}_{1}$:
$$\mathcal{T}_{1}\subset C^{*}(T_{\Phi_1},T_{\Phi_{\frac{m}{n}}})\cong\mathcal{T}_{\frac{1}{n}}.$$

Similarly, as was done in Lemma it can be shown:
$$\mathcal{T}_{\frac{1}{n}}\subset C^{*}(T_{\Phi_{\frac{1}{n}}},T_{\Phi_{\frac{m}{n^2}}})\cong\mathcal{T}_{\frac{1}{n^2}},$$
and for any $k$
$$\mathcal{T}_{\frac{1}{n^k}}\subset C^{*}(T_{\Phi_{\frac{1}{n^k}}},T_{\Phi_{\frac{m}{n^{k+1}}}})\cong\mathcal{T}_{\frac{1}{n^{k+1}}}.$$

Thus, we obtain the sequence of the Toeplitz $C^{*}$-algebras
$$\mathcal{T}_{1}\subset\mathcal{T}_{\frac{1}{n}}\subset\mathcal{T}_{\frac{1}{n^2}}\subset\mathcal{T}_{\frac{1}{n^3}}\subset\dots,$$
in which every next $C^{*}$-algebra is an inverse extension of the previous algebra.

\begin{theorem}
The inductive limit of Toeplitz $C^{*}$-algebras
$$\mathcal{T}_{1}\xrightarrow{j}\mathcal{T}_{\frac{1}{n}}\xrightarrow{j}\mathcal{T}_{\frac{1}{n^2}}\dots,$$
where $j$ is embedding, generates the $C^{*}$-algebra which is isomorphic to $C^{*}_{red}(\mathbb{Q}^{(n)}_{+})$, where $\mathbb{Q}^{(n)}$ is the semigroup of rational numbers generated by the numbers of the form $\frac{m}{n^k}$, where $m\in\mathbb{Z}, k\in\mathbb{N}$, and $\mathbb{Q}^{(n)}_{+}=\mathbb{Q}_{+}\cap \mathbb{Q}^{(n)}$.
\end{theorem}

Now consider the case when $t$ is irrational positive number. Let $\Gamma$ be a group generated by the numbers $m+nt$ which is everywhere dense in $\mathbb{R}$, where $m,n\in\mathbb{Z}$. Denote $\Gamma_{+}=\Gamma\cap\mathbb{R}_{+}$.

\begin{theorem}
Let $C^{*}(T_{\Phi_1},T_{\Phi_t})$ is the $C^{*}$-algebra generated by $T_{\Phi_{1}}$ and $T_{\Phi_{t}}$. Then $C^{*}(T_{\Phi_1},T_{\Phi_t})$ and $C^{*}_{red}(\Gamma_{+})$ are canonically isomorphic.
\begin{proof}
Suppose we have the representation ${\pi}:\Gamma_{+}\rightarrow B(H^2)$, given by the following way:

\[
m+nt \mapsto
\begin{cases}
T^{m}_{\Phi_1}T^{n}_{\Phi_t}, \ \mbox{ if } n>0,m>0; \\
(T^{*}_{\Phi_1})^{|m|}T^{n}_{\Phi_t}, \ \mbox{ if } n>0,m<0; \\
(T^{*}_{\Phi_t})^{|n|}T^{m}_{\Phi_1}, \ \mbox{ if } n<0,m>0.
\end{cases}
\]

Similarly to the theorem~\ref{tirat} it can be shown that
$$\pi(m+nt)=T_{\Phi_{m+nt}}$$
and, obviously, $T_{\Phi_{m+nt}}$ is the isometric multiplicative operator of the multiplication on inner function $\Phi_{m+nt} \ (m+nt>0)$. Note that $C^{*}(T_{\Phi_1},T_{\Phi_t})$ is generated by  isometric representation $\pi$.

Thus, it follows from the theorem of Douglas \cite{Douglas}, that $C^{*}(T_{\Phi_1},T_{\Phi_t})$ and $C^{*}_{red}(\Gamma_{+})$ are canonically isomorphic.
\end{proof}
\end{theorem}


{\small \vspace{\baselineskip}\hrule \vspace{3pt}
\par
{\bf Grigoryan Tamara Anatolevna} -- {Kazan State Power Engineering University}
\par
E-mail: {\it tkhorkova@gmail.com} }
{\small \vspace{\baselineskip}\hrule \vspace{3pt}
\par
{\bf Lipacheva Ekaterina Vladimirovna} -- Kazan State Power Engineering University
\par
E-mail: {\it elipacheva@gmail.com} }
{\small \vspace{\baselineskip}\hrule \vspace{3pt}
\par
{\bf Tepoyan Vardan Hakobi} -- Kazan State Power Engineering University
\par
E-mail: {\it tepoyan.math@gmail.com} }

\end{document}